\documentclass[10pt,a4paper,fleqn,draft]{article}
\usepackage{style,fullpage}

\let\*\bullet
\let\ot\otimes
\let\fn\rightarrow
\let\w\omega
\def\first{\hskip1cm&\hskip-1cm}
\def\mapf#1{\stackrel{#1}{\longmapsto}}
\newcommand{\cf}{\textit{cf.}}
\DeclareMathOperator{\Imagen}{Im}
\DeclareMathOperator{\Ker}{Ker}
\newcommand{\B}{\mathscr{B}}
\newcommand{\Sym}{\mathscr{S}}

\title{Twisted K\"ahler differential forms}
\author{%
Max Karoubi\thanks{Universit\'e Paris 7 - Math\'ematiques - 2, place
Jussieu 75251 Paris Cedex 05.}
\and
Mariano Suarez Alvarez\thanks{Departamento de Matem\'atica,
Facultad de Ciencias Exactas y Naturales, Universidad
de Buenos Aires. Ciudad Universitaria. Pab I, Buenos Aires (1428)
Argentina. e-mail: \texttt{mariano@dm.uba.ar}.\newline
\hspace*{\parindent} % why on earth isn't \thanks defined \long?
This work was supported by a grant from UBACYT TW69, the
international cooperation project SECyT-ECOS A98E05, and a CONICET
scholarship.}
}
\date{}

\begin{document}
\maketitle

%%%%%%%%%%%%%%%%%%%%%%%%%%%%%%%%%%%%%%%%%%%%%%%%%%
In \cite{Karoubi}, the first author has shown the interest of
``quantum'' differential forms in Algebraic Topology. They are obtained
from the usual ones by a slight change of the rules of calculus on
polynomials and series. In this paper, we make a more systematic 
study of these new
quantum differential forms. Our starting point is a commutative algebra $A$
with an endomorphism $\alpha$; the dif\-ferential graded algebra of
``twisted'' differential forms $\Omega^\*_\alpha A$ is then obtained as
the quotient of the universal non-commutative differential graded algebra
$\Omega^\* A$, defined by A.~Connes and the first author, by the ideal
generated by the relations ($d$ being the differential)
   \[
   da\,b-\alpha(b)\,da.
   \]
If $\alpha$ is the identity, we recover the classical commutative
differential graded algebra of K\"ahler differential forms. If $A=k[x]$ and
the endomorphism $\alpha$ is given by $\alpha(x^n)=q^nx^n$, where $q\in k$
is a ``quantum'' parameter, we find the differential graded algebra
introduced in \cite{Karoubi} for topological purposes.

The interest of this general definition lies essentially in the
existence of a remarkable braided structure $R$ on $\Omega^\*_\alpha A$,
which reduces to the ordinary flip if $\alpha$ is the identity, in the way
defined in \cite{Karoubi}, p.~2---see the precise definition below. As a
matter of fact, we show at the same time its \emph{uniqueness} under the
condition that $R(a\ot b)=b\ot a$ when both $a$ and $b$ belong to $A$,
identified to the degree zero part of $\Omega^\*_\alpha A$.  If $\alpha$ is
an automorphism, we produce in this way a lot of examples of representations
of the braid group $\B_n$ in a vector space or a module, by considering
$(\Omega^\*_\alpha A)^{\ot n}$ or, more generally, $J^{\ot n}$, where $J$
is any sub-quotient of $\Omega^\*_\alpha A$ stable by the braiding.  For
instance, if $A$ is the algebra of polynomials in several variables and if
$\alpha$ is induced by a linear transformation of these variables, 
filtrations by
various degrees in the variables produce such sub-quotients.

%%%%%%%%%%%%%%%%%%%%%%%%%%%%%%%%%%%%%%%%%%%%%%%%%%
\section{Generalities and statement of the theorem}

\paragraph Let $A$ be an associative algebra. A \emph{universal 
derivation} for $A$
is a derivation $d:A\fn\Omega^1A$ of $A$
such that for each derivation $\delta:A\fn M$ of $A$ with values in an
$A$-bimodule $M$ there exists exactly one morphism of bimodules
$f:\Omega^1A\fn M$ such that $\delta=f\circ d$. Such an object always
exists, and is unique up to an obvious notion of isomorphism; it can 
be concretely
realized by taking $\Omega^1A=\Ker(A\ot A\fn A)$, the kernel of the 
multiplication map,
and defining $da=1\ot a-a\ot1$ if $a\in A$.

\paragraph The \emph{algebra of universal differential forms} on $A$, 
which we shall  write
$\Omega^\*A$, is the tensor algebra
$T_A\Omega^1\!A$ of the $A$-bimodule $\Omega^1A$; it has a natural grading,
and the map $d:A\fn\Omega^1A$ induces in a unique way a derivation
$d:\Omega^\*A\fn\Omega^\*A$ with respect to which it becomes a
cohomologically graded differential algebra; \cf~\cite{Connes,Karoubi2}.

\paragraph Let now $A$ be a \emph{commutative} algebra, and let 
$\alpha:A\fn A$ be an algebra
endomorphism; we write $\Bar a=\alpha(a)$.
Let $I_\alpha A$ be the differential ideal generated in
$\Omega^\*A$ by the elements $da\,b-\Bar b\, da$ for $a$, $b\in A$,
and $\Omega_\alpha^\*A=\Omega^\*A/I_\alpha A$. This is again by construction a
differential  algebra, which is graded since the  ideal $I_\alpha$ is
homogeneous, and which is clearly natural with respect to maps in the
category of pairs $(A,\alpha)$ as above, and where the morphisms are 
morphisms of the
underlying algebras commuting with the given endomorphisms. We call
$\Omega_\alpha^\*A$ the \emph{differential graded algebra of twisted K\"ahler
differential forms}.

We note that since $I_\alpha$ is a differential ideal we have the relation
$dudv=-d\Bar vdu$ in $\Omega^\*_\alpha $
for each pair of elements $u,v\in A$, as a simple computation shows.

\paragraph Let $A$ be an algebra. A \emph{braiding} on $A$ is a 
morphism $R:A\ot
A\fn A\ot A$ such that
   \begin{gather}
    R\circ(\eta\ot 1)=1\ot\eta, \qquad R\circ(1\ot\eta)=\eta\ot1; 
\label{unit} \\
    (R\ot1)\circ(1\ot R)\circ (R\ot 1)=(1\ot R)\circ (R\ot1)\circ(1\ot 
R); \label{braid_eq} \\
    (\mu\ot1) \circ (1\ot R)\circ (R\ot 1)=R\circ(1\ot\mu); \label{compat:1} \\
    (1\ot\mu)\circ (R\ot 1)\circ (1\ot R)=R\circ(\mu\ot1); \label{compat:2} \\
    \mu\circ R=\mu.\label{mult}
   \end{gather}
Here $\mu:A\ot A\fn A$ is the multiplication map, and $\eta:k\fn A$ gives
the identity element of $A$.

The operator $R$ is regarded as an interchange operator. From this point of
view, the condition \eqref{braid_eq}, the Yang-Baxter equation, is a
natural one to impose; in particular, it implies that there is an action of
the braid group $\B_n$ on the tensor power $A^{\ot n}$ whenever $\alpha$ is
an automorphism.  Relations \eqref{compat:1} and \eqref{compat:2} express
compatibility of the braiding with the product. Finally, equation
\eqref{mult} is read as imposing a commutativity.

\paragraph If $A$ is a differential graded  algebra, we will say that
a morphism $R:A\ot A\fn A\ot A$ is a \emph{braiding of differential graded
algebras} if it is simultaneously a braiding and a map of
differential graded modules with respect to the usual structure on $A\ot A$.

\paragraph If $A$ is a commutative algebra, we consider the morphism 
$\tau:A\ot A\fn
A\ot A$ given by $\tau(a\ot b)=b\ot a$; it is a braiding of $A$, the
\emph{trivial braiding} or \emph{ordinary flip}.

\paragraph With this vocabulary, we can now state our theorem:

\begin{Theorem}\label{the:theorem}
There exists a unique functorial way of assigning to each 
endomorphism $\alpha:A\fn A$
of a commutative algebra a braiding
$R:\Omega_\alpha^\*A\ot\Omega_\alpha^\*A\fn\Omega_\alpha^\*A\ot\Omega_\alpha^\*A$
of the differential graded algebra of twisted K\"ahler differential forms
on $(A,\alpha)$ in such a way that its restriction to the degree zero
submodule $\Omega_\alpha^0A\ot\Omega_\alpha^0A=A\ot A$ is the trivial 
braiding $\tau$.
\end{Theorem}

%%%%%%%%%%%%%%%%%%%%%%%%%%%%%%%%%%%%%%%%%%%%%%%%%%
\section{Uniqueness}

\paragraph We write $R_{i,j}$ for the restriction of $R$ to $\Omega^i_\alpha
A\ot\Omega^j_\alpha A$. Our strategy to show uniqueness is to relate the
various restrictions $R_{i,j}$ to $R_{0,0}$ and $R_{1,0}$ and then to prove
that these two morphisms are determined by the conditions stated in 
the theorem.

We start with a straightforward computation.

\paragraph Let $i\geq0$, $j\geq1$. One has
   \begin{align*}\first
   R_{i,j}(u_0du_1\cdots du_i\ot v_0dv_1\cdots dv_j)\\
     =\;& -R(u_0du_1\cdots du_i\ot d(v_0v_1)\cdots dv_j)\\
      & -\sum_{k=1}^{j-1}(-1)^kR
          (u_0du_1\cdots du_i\ot dv_0dv_1\cdots d(v_kv_{k+1})\cdots dv_j)\\
      & -(-1)^jR(u_0du_1\cdots du_i\ot dv_0dv_1\cdots dv_{j-1}v_j)\\
     =\;& -(-1)^iRd (u_0du_1\cdots du_i\ot v_0v_1dv_2\cdots dv_j)
        +(-1)^iR(du_0du_1\cdots du_i\ot v_0v_1dv_2\cdots dv_j)\\
      & -\sum_{k=1}^{j-1}(-1)^{k+i}Rd
          (u_0du_1\cdots du_i\ot v_0dv_1\cdots d(v_kv_{k+1})\cdots dv_j)\\
      & +\sum_{k=1}^{j-1}(-1)^{k+i}R
          (du_0du_1\cdots du_i\ot v_0dv_1\cdots d(v_kv_{k+1})\cdots dv_j)\\
      & -(-1)^jR(u_0du_1\cdots du_i\ot dv_0dv_1\cdots dv_{j-1}v_j)\\
\intertext{}
     =\;& -(-1)^idR(u_0du_1\cdots du_i\ot
     	 v_0[v_1dv_2\cdots dv_j
	    +\sum_{k=1}^{j-1}(-1)^k)dv_1\cdots d(v_kv_{k+1})\cdots dv_j])\\
      & +(-1)^iR(du_0du_1\cdots du_i\ot
     	 v_0[v_1dv_2\cdots dv_j
	    +\sum_{k=1}^{j-1}(-1)^k)dv_1\cdots d(v_kv_{k+1})\cdots dv_j])\\
      & -(-1)^jR(u_0du_1\cdots du_i\ot dv_0dv_1\cdots dv_{j-1}v_j)\\
     =\;&  (-1)^{i+j}dR_{i,j-1}(u_0du_1\cdots du_i\ot v_0dv_1\cdots 
dv_{j-1}v_j)\\
      & -(-1)^{i+1}R_{i+1,j-1}(du_0du_1\cdots du_i\ot v_0dv_1\cdots 
dv_{j-1}v_j)\\
      & -(-1)^jR(u_0du_1\cdots du_i\ot dv_0dv_1\cdots dv_{j-1}v_j).
   \end{align*}
Since
   \begin{align*}\first
   R(u_0du_1\cdots du_i\ot dv_0\cdots dv_{j-1}v_j) \\
    =\;&(R\circ1\ot\mu)(u_0du_1\cdots du_i\ot dv_0\cdots dv_{j-1}\ot v_j) \\
    =\;&(\mu\ot1\circ1\ot R\circ R\ot1)
         (u_0du_1\cdots du_i\ot dv_0\cdots dv_{j-1}\ot v_j)\\
    =\;&(\mu\ot1\circ1\ot R_{*,0})
         (R(u_0du_1\cdots du_i\ot dv_0\cdots dv_{j-1})\ot v_j)\\
\intertext{and}
   \first R(u_0du_1\cdots du_i\ot dv_0\cdots dv_{j-1})\\
     =\;& (-1)^iRd(u_0du_1\cdots du_i\ot dv_0\cdots dv_{j-1})
       -(-1)^iR(du_0du_1\cdots du_i\ot v_0dv_1\cdots dv_{j-1})\\
     =\;& (-1)^idR_{i,j-1}(u_0du_1\cdots du_i\ot dv_0\cdots dv_{j-1})
       -(-1)^iR_{i+1,j-1}(du_0du_1\cdots du_i\ot v_0dv_1\cdots dv_{j-1})
   \end{align*}
we see that we can compute $R_{i,j}$ in terms of $R_{i,j-1}$, $R_{i+1,j-1}$ and
$R_{\*,0}$; a simple inductive argument then shows that $R$
is determined by $R_{\*,0}$.

\paragraph Now, if $i\geq0$,
   \begin{align*}
   \first R_{i+1,0}(u_0du_1\cdots du_{i+1}\ot v_0)
     = (R\circ\mu\ot1)(u_0du_1\cdots du_i\ot du_{i+1}\ot v_0)\\
     =\;&(1\ot\mu\circ R\ot1\circ1\ot R)(u_0du_1\cdots du_i\ot 
du_{i+1}\ot v_0)\\
     =\;&(1\ot\mu\circ R\ot1)(u_0du_1\cdots du_i\ot R_{1,0}(du_{i+1}\ot v_0)),
   \end{align*}
so that, if we assume
   \begin{equation}\label{imagen}
   \Imagen R_{1,0}\subset \Omega_\alpha^0A\ot\Omega_\alpha^1A,
   \end{equation}
we see that the maps $R_{\*,0}$ are determined by $R_{0,0}$ and $R_{1,0}$.

\paragraph Part of our hypothesis is that $R_{0,0}=\tau$; the required
uniqueness will follow then if we can show that the hypothesis also 
determines $R_{1,0}$ in
a such a way that \eqref{imagen} is verified.

\paragraph Let us consider the polynomial algebra 
$L_2=k[\{x_i,y_i\}_{i\geq0}]$ on
variables $x_i$ and $y_i$, for $i\geq 0$, equipped with the endomorphism
$\lambda:L_2\fn L_2$ such that $\lambda(x_i)=x_{i+1}$ and 
$\lambda(y_i)=y_{i+1}$.

Since $R$ is a braiding, we have that
$\mu R(dx_0\ot y_0) = \mu(dx_0\ot y_0) = dx_0y_0 = y_1dx_0 =
\mu(y_1\ot dx_0)$, so there is an element
$\w\in\Ker(\mu:\Omega_\lambda^\*L_2\ot\Omega_\lambda^\*L_2\fn\Omega_\lambda^\*L_2)$
such that
   \[
   R(dx_0\ot y_0)=y_1\ot dx_0+\w.
   \]

\paragraph Let $A$ be a commutative algebra and $\alpha:A\fn A$ be an 
endomorphism of
$A$; if $a$,~$b\in A$, there is exactly one morphism in the category of
endomorphisms of algebras $\phi_{a,b}:(L_2,\lambda)\fn(A,\alpha)$
such that $\phi_{a,b}(x_0)=a$ and $\phi_{a,b}(y_0)=b$,
and it induces in turn a morphism of differential graded algebras, 
which we will
write $\phi_{a,b}$ as well, 
$\phi_{a,b}:\Omega_\lambda^\*L_2\fn\Omega_\alpha^\*A$.
Naturality of $R$ implies that
   \begin{align*}\first
   R(da\ot b) = R\phi_{a,b}(dx_0\ot y_0) =\phi_{a,b}R(dx_0\ot y_0) \\
     =\;&\phi_{a,b}(y_1\ot dx_0)+\phi_{a,b}(\w) \\
     =\;&\Bar b\ot da+\phi_{a,b}(\w);
   \end{align*}
so that $\w$ determines $R_{1,0}$ on the elements of the form $da\ot b$
in $\Omega_\alpha^1A\ot\Omega_\alpha^0A$. In general, if
$adb\ot c\in \Omega_\alpha^1A\ot\Omega_\alpha^0A$, we have
   \begin{align*}\first
   R(adb\ot c) = R(d(ab)\ot c)-R(da b\ot c) \\
     =\;& R(d(ab)\ot c)-(R\circ\mu\ot1)(da\ot b\ot c) \\
     =\;& R(d(ab)\ot c)-(1\ot\mu\circ R\ot1\circ1\ot R)(da\ot b\ot c) \\
     =\;& R(d(ab)\ot c)-(1\ot\mu\circ R\ot1)(da\ot c\ot b) \\
     =\;& R(d(ab)\ot c)-(1\ot\mu)(R(da\ot c)\ot b);
   \end{align*}
observe that we have used the hypothesis that $R_{0,0}=\tau$.
We conclude that $\w$ actually determines $R_{1,0}$.

Let us write $\w(a,b)=\phi_{a,b}(\w)$.

\paragraph Now let $L_3=k[\{x_i,y_i,z_i\}_{i\geq0}]$
be endowed with the endomorphism $\lambda:L_3\fn L_3$
such that $\lambda(x_i)=x_{i+1}$, $\lambda(y_i)=y_{i+1}$ and 
$\lambda(z_i)=z_{i+1}$.
We compute in $\Omega_\lambda^\*L_3$:
   \begin{align*}\first
   (1\ot R\circ R\ot1\circ1\ot R)(dx_0\ot y_0\ot z_0)\\
   =\;& (1\ot R\circ R\ot 1)(dx_0\ot z_0\ot y_0)\\
   =\;& (1\ot R)(z_1\ot dx_0\ot y_0+\w(x_0,z_0)\ot y_0)\\
   =\;& z_1\ot y_1\ot dx_0+z_1\ot\w(x_0,y_0)+(1\ot R)(\w(x_0,z_0)\ot y_0)\\
   \first
   (R\ot1\circ 1\ot R\circ R\ot1)(dx_0\ot y_0\ot z_0)\\
   =\;& (R\ot1\circ1\ot R)(y_1\ot dx_0\ot z_0+\w(x_0,y_0)\ot z_0)\\
   =\;&(R\ot1)(y_1\ot z_1\ot dx_0+y_1\ot\w(x_0,z_0)+(1\ot R)
   	(\w(x_0,y_0)\ot z_0)\\
   =\;& z_1\ot y_1\ot dx_0+(R\ot 1)(y_1\ot\w(x_0,z_0))
         + (R\ot1\circ1\ot R)(\w(x_0,y_0)\ot z_0)
   \end{align*}
Since $R$ satisfies the braid equation \eqref{braid_eq}, we have then that
   \[
   z_1\ot\w(x_0,y_0)+(1\ot R)(\w(x_0,z_0)\ot y_0)
   = (R\ot 1)(y_1\ot\w(x_0,z_0)) + (R\ot1\circ1\ot R)(\w(x_0,y_0)\ot z_0)
   \]
Apply $1\ot\mu$ to both sides of this equality; on the left, we obtain
   \begin{align*}\first
   (1\ot\mu)(z_0\ot\w(x_0,y_0))+(1\ot\mu R)(\w(x_0,z_0)\ot y_0) \\
   =\;&(1\ot\mu)(\w(x_0,z_0)\ot y_0)\\
   =\;&\w(x_0,z_0)y_0\\
\intertext{and, on the right,}
   \first
   (1\ot\mu \circ R\ot 1)(y_1\ot\w(x_0,z_0))
     + (1\ot\mu\circ R\ot1\circ1\ot R)(\w(x_0,y_0)\ot z_0) \\
   =\;& (1\ot\mu\circ R\ot1)(y_1\ot\w(x_0,z_0))
     + (\mu\ot1)(\w(x_0,y_0)\ot z_0) \\
   =\;& (1\ot\mu\circ R\ot 1)(y_1\ot\w(x_0,z_0)).
   \end{align*}
so that
   \[
   \w(x_0,z_0)y_0=(1\ot\mu\circ R\ot1)(y_1\ot \w(x_0,z_0)).
   \]
Observe that the variable $y_0$ cannot appear on the right hand side of this
equality because of naturality; in view of the left hand side, we must have
$\w=0$.

This shows that, if $\alpha:A\fn A$ is an endomorphism of a commutative
algebra, we have in $\Omega_\alpha^\*A$ that
   \[
   R(da\ot b)=\Bar b\ot da.
   \]
In view of what has been said above, the uniqueness statement in the
theorem follows from this.

%%%%%%%%%%%%%%%%%%%%%%%%%%%%%%%%%%%%%%%%%%%%%%%%%%
\section{Existence}

\paragraph Let us show now that there exists a braiding satisfying 
the conditions in
the statement. We do this by explicitly constructing it.

\paragraph Let $A$ be a commutative and let $\alpha:A\fn A$ be a 
endomorphism of $A$.
We define a morphism of graded modules
$I:\Omega_\alpha^\*A\fn\Omega_\alpha^\*A$
of degree $-1$ by putting, on $\Omega_\alpha^nA$,
   \[
   I(u_0du_1\cdots du_n)=
     \sum_{i=1}^n(-1)^{i+1} u_0 du_1 \cdots du_{i-1}(u_i-\Bar u_i)d\Bar u_{i+1}
     	\cdots d\Bar u_n.
   \]
It is easy to check that this is well defined. This operator is obviously
left $A$-linear, and does not
commute in general with the differential on $\Omega_\alpha^\*A$; in fact,
   \begin{align*}\first
   Id(u_0du_1\cdots du_n)
     = \sum_{i=0}^n(-1)^{i+2}du_0\cdots du_{i-1}(u_i-\Bar u_i)d\Bar 
u_{i+1}\cdots d\Bar u_n \\
\intertext{and}
   \first dI(u_0du_1\cdots du_n)
     = \sum_{i=1}^n(-1)^{i+1}d(u_0 \cdots du_{i-1}(u_i-\Bar u_i)d\Bar 
u_{i+1} \cdots
     		d\Bar u_n) \\
     =\;&\sum_{i=1}^n(-1)^{i+1}du_0 \cdots du_{i-1}(u_i-\Bar u_i)d\Bar 
u_{i+1} \cdots
                 d\Bar u_n \\
        &+\sum_{i=1}^n(-1)^{i+1}u_0 \cdots du_{i-1}(du_i-d\Bar 
u_i)d\Bar u_{i+1} \cdots
                 d\Bar u_n \\
\intertext{so that}
   \first (Id+dI)(u_0du_1\cdots du_n)
     = (u_0-\Bar u_0)d\Bar u_1\cdots d\Bar u_n
       +\sum_{i=1}^nu_0du_1\cdots du_{i-1}(du_i-d\Bar u_i)d\Bar 
u_{i+1}\cdots d\Bar u_n\\
     =\;& -\Bar u_0d\Bar u_1\cdots d\Bar u_n+u_0du_1\cdots du_n \\
     =\;& (1-\alpha)(u_0du_1\cdots du_n)
   \end{align*}
Thus, $I$ is a homotopy $1_{\Omega_\alpha^\*A}\simeq \alpha$.

\paragraph This computation proves the first part of the following 
lemma. To state it
and in order to simplify future formulas, we introduce some notation.  In
what follows we shall write $[n]$, for $n\in\mathbb Z$, instead of $(-1)^n$,
and, in a context where there is an endomorphism of an algebra---$\alpha$,
say---we will write $\Bar a^n$ instead of $\alpha^n(a)$. Also, we will
agree that a homogeneous differential form stands for its degree when
inside square brackets or in an exponent. For example, we will write
$[\w(\psi+1)]\Bar\phi^\w$, when $\w$,~$\psi$~and~$\phi$ are homogeneous
differential forms, instead of
$(-1)^{\deg\w(\deg\psi+1)}\alpha^{\deg\w}(\phi)$.

\begin{Lemma}\label{the:lemma}
Let $A$ be a commutative algebra and let $\alpha:A\fn A$ be an
endomorphism of $A$.
For $\w,\psi\in\Omega_\alpha^\*A$ and $v\in A$ the following relations
hold:
   \begin{gather}
   I:1_{\Omega_\alpha^\*A}\simeq \alpha\label{lema:i}\\
   I(\w\psi)=I\w\Bar\psi+[\w]\w I\psi\label{lema:ii}\\
   I(\w dv)=I\w d\Bar v+[\w]\w(v-\Bar v)\label{lema:iii}\\
   \w v-v\w=[\w]I\w dv\label{lema:iv}
   \end{gather}
Moreover, we have $I^2=0$.
\end{Lemma}

\begin{Proof}
We have just shown \eqref{lema:i}; \eqref{lema:ii} and \eqref{lema:iii}
follow immediately from the definition of $I$. Let us check \eqref{lema:iv}
inductively on $\deg\w$. If $\deg\w=0$, it is true because $A$ is
commutative and $I$ is zero on $0$-forms. Assume then the truth of
\eqref{lema:iv} for an homogeneous form $\w$; then, if $u$,~$v\in A$,
   \begin{align*}\first
   \w du v
     = \w d(uv)-\w udv \\
     =\;& \w vdu + \w dv u -\w u dv \\
     =\;& \w v du - \w (u - \Bar u) dv \\
     =\;& v \w du + [\w] I\w dv -\w(u-\Bar u) dv \\
     =\;& v \w du + [\w] (I\w-[\w]\w(u-\Bar u))dv \\
     =\;& v \w du + [\w] I(\w du) dv
   \end{align*}
This shows \eqref{lema:iv} for all $\w$.

Finally, to show that $I^2=0$ inductively, we observe that it is
trivially true on $0$-forms, and if $I^2\w=0$ for an homogeneous form $\w$,
we have
   \begin{align*}\first
   I^2(\w dv)
     = I(I\w d\Bar v+[\w]\w(v-\Bar v)) \\
     =\;& I^2\w d\Bar v^2 + [I\w]I\w I(d\Bar v)
          + [\w] I\w(\Bar v-\Bar v^2)+[2\w]\w I(v-\Bar v) \\
     =\;& [\w-1]I\w(\Bar v-\Bar v^2)+[w]I\w(\Bar v-\Bar v^2) \\
     =\;& 0
   \end{align*}
so that $I^2$ vanishes identically on $\Omega_\alpha^\*A$.~\qed
\end{Proof}

\paragraph Let us fix a commutative algebra $A$ and an endomorphism 
$\alpha:A\fn A$. Let
$R:\Omega_\alpha^\*A\ot\Omega_\alpha^\*A\fn\Omega_\alpha^\*A\ot\Omega_\alpha^\*A$
be given by
   \[
   R(\w\ot\phi)=[\w\phi]\Bar\phi^\w\ot\w-
   		   [(\w+1)\phi]I\Bar\phi^\w\ot d\w
   \]
We will verify that this operator satisfies the conditions in the theorem.
 From the definition, it is clear that $R$ verifies \eqref{unit}.

\paragraph Next, we have
   \begin{align*}\first
   \mu R(u_0du_1\cdots du_n\ot v_0dv_1\cdots dv_m) \\
     =\;&[nm]\Bar v_0^nd\Bar v_1^n\cdots d\Bar v_m^nu_0du_1\cdots du_n
        -[(n+1)m]I(\Bar v_0^nd\Bar v_1^n\cdots d\Bar v_m^n)
          du_0du_1\cdots du_n \\
     =\;&[nm]\Bar v_0^nu_0d\Bar v_1^n\cdots d\Bar v_m^ndu_1\cdots du_n
        +[(n+1)m]\Bar v_0^nI(d\Bar v_1^n\cdots d\Bar v_m^n)
          du_0du_1\cdots du_n \\
        &
        -[(n+1)m]I(\Bar v_0^nd\Bar v_1^n\cdots d\Bar v_m^n)
          du_0du_1\cdots du_n \\
     =\;&u_0du_1\cdots du_nv_0dv_1\cdots dv_m
   \end{align*}
so that $\mu R=\mu$; this is \eqref{mult}.

\paragraph We want to check that $R$ satisfies the braid equation
\eqref{braid_eq}; evaluating both sides on $\w\ot\phi\ot\psi$ for
homogeneous forms $\w$, $\phi$, $\psi\in\Omega_\alpha^\*A$, we find
   \begin{align*}\first
   \w\ot\phi\ot\psi \\
     \mapf{R\ot1}\;&
        [\phi\w]\Bar\phi^\w\ot\w\ot\psi
        -[\phi(\w+1)]I\Bar\phi^\w\ot d\w\ot\psi \\
      \mapf{1\ot R}\;&
         [\phi\w+\w\psi]\Bar\phi^\w\ot\Bar\psi^\w\ot\w \\
       &-[\phi\w+\w\psi+\psi]\Bar\phi^\w\ot
          	I\Bar\psi^\w\ot d\w \\
       &-[\phi\w+\phi+(\w+1)\psi]
           I\Bar\phi^\w\ot\Bar\psi^{\w+1}\ot d\w \\
      \mapf{R\ot 1}\;&
         [\phi\w+\w\psi+\phi\psi]
          \Bar\psi^{\w+\phi}\ot\Bar\phi^\w\ot\w \\
       &-[\phi\w+\w\psi+\phi\psi+\psi]
	  I\Bar\psi^{\w+\phi}\ot d\Bar\phi^\w\ot\w \\
       &-[\phi\w+\w\psi+\psi+\phi(\psi-1)]
	  I\Bar\psi^{\w+\phi}\ot\Bar\phi^\w\ot d\w \\
       &-[\phi\w+\phi+(\w+1)\psi+(\phi-1)\psi]
           \Bar\psi^{\w+\phi}\ot I\Bar\phi^\w\ot d\w \\
       &+[\phi\w+\phi+(\w+1)\psi+(\phi-1)\psi+\psi]
           I\Bar\psi^{\w+\phi}\ot dI\Bar\psi^\w\ot d\w\\
\intertext{and}
   \first
   \w\ot\phi\ot\psi \\
     \mapf{1\ot R}\;&
       [\phi\psi]\w\ot\Bar\psi^\phi\ot\phi  \\
       &-[\phi\psi+\psi]\w\ot I\Bar\psi^\phi\ot d\phi \\
     \mapf{R\ot1}\;&
       [\phi\psi+\w\psi]\Bar\psi^{\phi+\w}\ot\w\ot\phi \\
       &-[\phi\psi+\psi+\w(\psi-1)]
         I\Bar\psi^{\phi+\w}\ot\w\ot d\phi \\
     \mapf{1\ot R}\;&
       [\phi\psi+\w\psi+\w\phi]
          \Bar\psi^{\phi+\w}\ot\Bar\phi^\w\ot\w \\
        &-[\phi\psi+\w\psi+\w\phi+\phi]
          \Bar\psi^{\phi+\w}\ot I\Bar\phi^\w\ot d\w \\
        &-[\phi\psi+\w\psi+\w\phi+\phi]
          I\Bar\psi^{\phi+\w}\ot\Bar\phi^{\w+1}\ot d\w \\
        &-[\phi\psi+\psi+(\w(\psi-1)+\w(\phi+1)]
          I\Bar\psi^{\phi+\w}\ot d\Bar\phi^\w\ot\w \\
        &-[\phi\psi+\psi+\w(\psi-1)+\w(\phi+1)+\phi+1]
          I\Bar\psi^{\phi+\w}\ot Id\Bar\phi^\w\ot d\w
   \end{align*}
These are equal, because
   \begin{multline*}
   -[\phi\w+\w\psi+\psi+\phi\psi-\phi]
     I\Bar\psi^{\w+\phi}\ot\Bar\phi^\w\ot d \w
   +[\phi\w+\phi+\w\psi+\psi+\phi\psi]
     I\Bar\psi^{\w+\phi}\ot dI\Bar\phi^\w\ot d\w
     = \\
   -[\phi\psi+\w\psi+\psi+\w\phi+\phi]
     I\Bar\psi^{\phi\w}\ot\Bar\phi^{\w+1}\ot d\w
   +[\phi\psi+\psi+\w\psi+\w\phi+\phi+1]
     I\Bar\psi^{\phi+\w}\ot Id\Bar\phi^\w\ot d\w
   \end{multline*}
which in turn follows from
   \[
   -[\psi-\phi]\Bar\phi^\w+[\phi+\psi]dI\Bar\phi^\w
   = -[\psi+\phi]\Bar\phi^{\w+1}
     -[\psi+\phi]Id\Bar\phi^\w
   \]
which is true, because lemma \pref{the:lemma} implies that
   \[
   -\Bar\phi^\w+dI\Bar\phi^\w
   =
   -\Bar\phi^{\w+1}-Id\Bar\phi^\w
   \]

\paragraph Finally, our map $R$  is compatible with multiplication in 
$\Omega_\alpha^\*A$,
since, for example,
   \begin{align*}\first
   \w\ot\phi\ot\psi \\
     \mapf{1\ot\mu}\;& \w\ot\phi\psi \\
     \mapf{R}\;&
        [\w(\phi+\psi)]
          \Bar\phi^\w\Bar\psi^\w\ot\w
       -[\w(\phi+\psi)+\phi+\psi]
          I(\Bar\phi^\w\Bar\psi^\w)\ot d\w \\
   \first \w\ot\phi\ot\psi \\
     \mapf{R\ot 1}\;&
       [\w\phi]\Bar\phi^\w\ot\w\ot\psi
       -[\w\phi+\phi]I\Bar\psi^\w\ot d\w\ot\psi \\
     \mapf{1\ot R}\;&
       [\w\phi+\w\psi]
         \Bar\phi^\w\ot\Bar\psi^\w\ot\w \\
       &-[\w\phi+\w\psi+\psi]
         \Bar\phi^\w\ot I\Bar\psi^\w\ot d\w \\
       &-[\w\phi+\phi+(\w+1)\psi]
         I\Bar\phi^\w\ot\Bar\psi^{\w+1}\ot d\w \\
      \mapf{\mu\ot1}\;&
        [\w\phi+\w\psi]\Bar\phi^\w\Bar\psi^\w\ot\w \\
        &-[\w\phi+\w\psi+\psi+\psi]
         ([\phi]\Bar\phi^\w I\Bar\psi^\w
	 +I\Bar\phi^\w\Bar\psi^{\w+1})\ot d\w
   \end{align*}
and these are equal by the lemma; this shows \eqref{compat:1}, and the
other equation \eqref{compat:2} is checked in the same way.

\paragraph It is obvious that $R$ depends naturally on $\alpha$, and
reduces to the trivial twist $\tau$ in degree $0$. Since it satisfies the
required conditions, theorem \pref{the:theorem} is proved.

\paragraph \textbf{Two simple examples.} Consider $A=k[x]$ and $q\in k$,
and let $\alpha:A\fn A$ be the endomorphism such that $\alpha(x)=qx$. Then we
have $\Omega^0_\alpha A=k[x]$, $\Omega^1_\alpha A=k[x]dx$ and the twisted
exterior differential is given by $dx^n=n_q x^{n-1}dx$ for each $n\geq1$,
where, for each $n$, we define the $q$-integer $n_q=(1-q^n)/(1-q)$ when
$q\neq 1$, and $n_1=n$. The
operator $I:\Omega^1_\alpha A\fn\Omega^0_\alpha A$ is given by
$q$-integration of forms: $I(x^ndx)=(1-q)x^{n+1}$. Using this, we easily
obtain the following formulas for the braiding constructed above on
$\Omega^\$_\alpha A$:
   \begin{gather*}
   R(x^n\ot x^m) = x^m\ot x^n\\
   R(x^ndx\ot x^m) = q^mx^m\ot x^ndx\\
   R(x^n\ot x^mdx) = x^mdx\ot x^n+(1-q)x^{m+1}\ot x^{n-1}dx\\
   R(x^ndx\ot x^mdx) = -q^{m+1}x^mdx\ot x^ndx
   \end{gather*}
We thus recover the main example considered in \cite{Karoubi}. More
generally, one can replace $A$ with a ring of convergent power series $f$
with the endomorphism given by $\alpha(f)(x)=f(qx)$ like in
\cite{Karoubi3}.

\paragraph Another familiar example is the following. Let $A=k[x]/(x^2-x)$
and let $\alpha:A\fn A$ be such that $\alpha(x)=1-x$. Then
$\Omega_\alpha^\*A=\Omega^\*A$ can be identified with the differential
graded algebra of normalized cochains on the simplicial set $\{0,1\}$.
Since $\alpha^2=1$, the action of the braid group $\B_n$ reduces to the
action of the symmetric group $\Sym_n$.

%%%%%%%%%%%%%%%%%%%%%%%%%%%%%%%%%%%%%%%%%%%%%%%%%%%%%%%%%%%%%%%%%%%%%%

\end{document}